# Convergence in Energy-Lowering (Disordered) Stochastic Spin Systems

**Emilio De Santis**

desantis@mat.uniroma1.it

*Università di Roma La Sapienza, Dipartimento di Matematica "Guido Castelnuovo"*

*Piazzale Aldo Moro, 2 - 00185 Roma, Italia*

**Charles M. Newman**

newman@courant.nyu.edu

*Courant Inst. of Mathematical Sciences, New York University, New York, NY 10012, USA*

October 31, 2018

**Abstract**

We consider stochastic processes, $S^t \equiv (S^t_x : x \in \mathbb{Z}^d) \in \mathcal{S}_0^{\mathbb{Z}^d}$ with $\mathcal{S}_0$ finite, in which spin flips (i.e., changes of $S^t_x$) do not raise the energy. We extend earlier results of Nanda-Newman-Stein that each site $x$ has almost surely only finitely many flips that strictly lower the energy and thus that in models without zero-energy flips there is convergence to an absorbing state. In particular, the assumption of finite mean energy density can be eliminated by constructing a percolation-theoretic Lyapunov function density as a substitute for the mean energy density. Our results apply to random energy functions with a translation-invariant distribution and to quite general (not necessarily Markovian) dynamics.

## 1 Introduction and Main results

In the statistical physics of Ising models, Glauber dynamics refers to Markov processes whose transitions are the flips of single spin variables and are such that the Gibbs distribution at a given nonzero temperature is invariant. The zero-temperature limit of these processes are interacting particle systems in which each transition lowers or leaves unchanged the total energy of the system. In [9], Nanda, Newman and Stein studied these zero-temperature processes for (homogeneous and) disordered Ising models on $\mathbb{Z}^d$ (and other lattices) with Hamiltonian (i.e., total energy)

$$\mathcal{H}(S) = \sum_{x \in \mathbb{Z}^d} V_x(\{S_z : ||z - x|| \le 1\}), \qquad (1)$$

where

$$V_x = -\frac{1}{2} \sum_{z:\,||z-x||=1} J_{\{x,z\}} S_x S_z. \tag{2}$$

Here $||\cdot||$ denotes Euclidean norm, the $S_x$'s are $\pm 1$ valued and the $J_{\{x,z\}}$'s are i.i.d. random variables (that are fixed as the dynamics runs its course). The initial spin configuration $S^0 \equiv (S_x^0 : x \in \mathbb{Z}^d)$ was chosen from (i.i.d.) product measure (with $P(S_x^0 = +1) = \lambda$) and the continuous-time Markov process $S^t$ has each site $x$ flip with rate either $1$, $1/2$ or $0$ depending on whether the energy change $\Delta_x \mathcal{H}$ caused by the flip is $< 0$, $= 0$, or $> 0$. We remark that throughout this paper sums over $\mathbb{Z}^d$ as in (1) are formal expressions while quantities such as $\Delta_x \mathcal{H}$ are well-defined and given by finite sums.

Under the assumption that

$$E|J_{\{x,z\}}| < \infty, \tag{3}$$

it was proved in [9] that a.s. each site $x$ has only finitely many flips with $\Delta_x \mathcal{H} < 0$. As a corollary, it follows that if the common distribution of the $J_{\{x,z\}}$'s is continuous (so that there is zero probability for the $J_{\{x,z\}}$'s to provide for a flip with $\Delta_x \mathcal{H} = 0$), then a.s. each site flips only finitely many times and so a.s. $S^t$ converges (coordinatewise) to some $S^\infty$ (an absorbing or metastable state). It was also shown in [9] that the finite mean energy density assumption (3) was not needed for $d = 1$ and could *sometimes* be avoided for $d \geq 2$. Our first main result is the next theorem which eliminates the assumption (3) in any dimension. The proof, in which the mean energy density is replaced by a percolation-theoretic Lyapunov function, is perhaps as interesting as the theorem itself. We slightly generalize the form of (2) since it is basically cost-free.

**Theorem 1.** *Let $S^t$ be the stochastic Ising model just described, except that (2) is generalized to*

$$V_x = -h_x S_x - \frac{1}{2} \sum_{z:\,||z-x||=1} J_{\{x,z\}} S_x S_z, \tag{4}$$

*with the $J_{\{x,z\}}$'s and $h_x$'s i.i.d. (and independent of each other and with common distributions $F_1$ and $F_2$). For any $F_1$, $F_2$ and any dimension $d$, a.s. each site $x$ has only finitely many flips with $\Delta_x \mathcal{H} \neq 0$. If either $F_1$ or $F_2$ is continuous, then a.s. each site has only finitely many flips of any kind.*

The last claim of Theorem 1 follows because the energy change $\Delta_x \mathcal{H}(S)$ caused by $S_x \to -S_x$ (coming from the changes in all the $V_y$'s with $||y - x|| \leq 1$) is

$$\Delta_x \mathcal{H} = 2h_x S_x + 2 \sum_{z:\,||z-x||=1} J_{\{x,z\}} S_x S_z, \tag{5}$$

and, for $F_1$ or $F_2$ continuous, this has zero probability of vanishing (for any of the $2^{2d+1}$ possible sign-assignments to $S_x$ and the $S_z$'s). The proof of the main part of Theorem 1 is given in Section 2 below as a corollary of a more general result, Theorem 2.

**Remark 1.1.** *As in [9], the proof and conclusions of Theorem 1 extend to models on much more general lattices and graphs than $\mathbb{Z}^d$. These include not only regular d-dimensional lattices (such as the two dimensional hexagonal lattice) but (as in [3]) also homogeneous trees and Cayley graphs of finitely generated groups.*

**Remark 1.2.** *There are many interesting cases where the final conclusion of Theorem 1 (that $S^t \to S^\infty$ a.s.) applies because $\Delta_x \mathcal{H}(S) \neq 0$ a.s. for any $S$ for a different reason than continuity of $F_1$ or $F_2$. For example, this is the case when $J_{\{x,z\}} \equiv J > 0$ and $h_x \equiv 0$ on graphs where every site has only an odd number of neighbors (such as the hexagonal lattice). It is also so on $\mathbb{Z}^d$ when $J_{\{x,z\}} \equiv J > 0$ and $h_x \equiv h$ provided $h/J \notin \{-2d, -2d+2, \ldots, 2d-2, 2d\}$.*

**Remark 1.3.** *There are also very interesting cases where $S^t \to S^\infty$ a.s. even though zero-energy flips do occur. This has been proved [6] (with $S^\infty \equiv 1$ or $\equiv -1$) on $\mathbb{Z}^d$ with $d \geq 2$ when $J_{\{x,z\}} \equiv J > 0$ and $h_x \equiv 0$ providing $\lambda \equiv P(S_x^0 = +1)$ is close enough to 1 or 0. When $\lambda = 1/2$, it has been proved for $d = 1$ and $d = 2$ (see, resp., [1] and [9]) that a.s. each site flips infinitely many times. Some other results about "local recurrence" on $\mathbb{Z}^2$ may be found in [4]*

**Remark 1.4.** *On $\mathbb{Z}^d$, with $h_x \equiv 0$ and $\lambda = 1/2$, the pair $(d, F_1)$ is said to be of type $\mathcal{I}$ (resp., $\mathcal{F}$) if a.s. each site flips infinitely (resp., only finitely) many times (see [7]). It was shown in [7] that $(d = 2, F_1 = \alpha \delta_J + (1-\alpha)\delta_{-J})$ with $J > 0$ and $0 < \alpha < 1$ is of (mixed) type $\mathcal{M}$, which means that a.s. some sites flip infinitely many times and some only finitely many. Combining all the previously known results (see [9, 7]) with Theorem 1, one finds that all $(d, F_1)$ have now been classified as type $\mathcal{I}$, $\mathcal{F}$ or $\mathcal{M}$ with the exception of the very interesting cases of $d \geq 3$ and $F_1 = \alpha \delta_J + (1-\alpha)\delta_{-J}$ with $J > 0$ and $0 \leq \alpha \leq 1$. These include the homogeneous ferromagnet ($\alpha = 1$) and antiferromagnet ($\alpha = 0$), which can easily be seen to be of the same type, and the $\pm J$ spin glass ($\alpha = 1/2$).*

Theorem 1 can immediately be extended in other useful ways besides those of Remark 1.1. For example, the initial distribution of $S^0$ can be any translation invariant measure and $h_x$ need not be independent of the $J_{\{x,z\}}$'s. Such cases will be covered by Theorem 2 below, which extends Theorem 1 in three basic but far-reaching ways:

(A) the spin variables $S_x$ need not be $\pm 1$ but can take their values from some arbitrary finite set $\mathcal{S}_0$;

(B) the stochastic process $(S^t : t \geq 0)$ (including the initial distribution for $S^0$) can be much more general than the specific zero-temperature Glauber dynamics Markov process with i.i.d. $S_x^0$'s treated in Theorem 1; e.g., it need not even be a Markov process;

(C) the random $V_x$'s of the Hamiltonian need not come from i.i.d. couplings but only need have a translation-invariant distribution.

We give a few natural examples of such extensions at the end of this section.

Needless to say, Theorem 2 can also be extended. One extension, to finite-range rather than only nearest-neighbor Hamiltonians, will be described in Remark 1.5 below. Other extensions, that we will not discuss, are to random lattices/graphs and to dynamics where multiple spin changes occur simultaneously (including continuous time Kawasaki dynamics and some discrete-time synchronous dynamics).

There are four main assumptions on the stochastic process $\hat{S} \equiv (S^t : t \geq 0)$ and random Hamiltonian $\mathcal{H}$ (determined by $(V_x : x \leq 0)$) which lead to the conclusion that a.s. each site $x$ flips (i.e. $S_x$ changes its value) with $\Delta_x \mathcal{H} \neq 0$ only finitely many times. The first three of these are fairly natural and unobjectionable: (1) statistical translational invariance, (2) single spin flip dynamics and (3) zero-temperature dynamics—i.e., $\Delta_x \mathcal{H} \leq 0$. The fourth assumption is that *either* of two hypothesis on $\mathcal{H}$ be valid: (4a) finite mean energy density or (4b) a somewhat technical looking large deviation bound on a (dependent) percolation model related to $\mathcal{H}$. This assumption, which is the price we pay for the generality of the theorem, did not appear in Theorem 1 because there the i.i.d. $J_{\{x,z\}}$'s result in a standard independent bond percolation model and the large deviation bound follows from standard results.

The finite mean energy density hypothesis 4a is the analogue of (3), the finite mean coupling hypothesis of [9], and the proof of Theorem 2 in that case is basically taken from [9]. On the other hand, the proof under the percolation hypothesis 4b involves a novel extension of the methods of [9], whereby a percolation-theoretic Lyapunov function, with finite mean density, is constructed as a replacement for the mean energy density. We state the percolation hypothesis and Theorem 2 in the context of nearest-neighbor Hamiltonians, where each $V_x$ is a function only of $\{S_z : ||z-x|| \leq 1\}$; this nearest-neighbor restriction, not needed at all with finite mean energy density, is included to simplify the statement and use of the percolation hypothesis (see Remark 1.5).

*The strong percolation hypothesis.* A slightly stronger than necessary (see Remark 1.6) version of the percolation hypothesis is as follows. Writing $V_y(S_x = \eta)$ to denote the function $V_y$ of $\{S_z : ||z - y|| \leq 1, z \neq x\}$ when $S_x$ is set to the value $\eta$, we define for $||x - y|| = 1$ the random variable

$$K_{x,y} = \max_{\{S_z : ||z-y|| \leq 1, z \neq x\}} \max_{\eta,\eta'} \{|V_y(S_x = \eta) - V_y(S_x = \eta')|\} \tag{6}$$

and

$$K^*_{\{x,y\}} = \max\{K_{x,y}, K_{y,x}\}. \tag{7}$$

For each fixed $K > 0$, we consider the (dependent) percolation model in which nearest-neighbor bonds $\{x, y\}$ with $||x - y|| = 1$ are said to be open if $K^*_{\{x,y\}} > K$ and otherwise closed and denote by $C_x^K$ the open cluster containing the site $x$ and by $|C_x^K|$ the number of sites in $C_x^K$. We will say that $\mathcal{H}$ satisfies the strong percolation hypothesis if there exist constants $M_K < \infty$ and $\lambda(K) \geq 0$ with $\lim_{K \to \infty} \lambda(K) = \infty$ such that for all $x \in \mathbb{Z}^d$

$$P(|C_x^K| > n) \leq M_K e^{-\lambda(K)n}. \tag{8}$$

In Theorem 2 we consider a stochastic process $\hat{S} \equiv (S^t : t \geq 0)$ taking values in $\mathcal{S}_0^{\mathbb{Z}^d}$, where $\mathcal{S}_0$ is finite, and a random Hamiltonian of the form

$$\mathcal{H}(S) = \sum_{x \in \mathbb{Z}^d} V_x(\{S_z : ||z - x|| \leq 1\}) \tag{9}$$

that satisfy the following three hypotheses:

1. *Translation invariance.* There is translation invariance for the joint distribution of $(\hat{S}, (V_x : x \in \mathbb{Z}^d))$.

2. *Single spin flip dynamics.* Almost surely, for each $x$, $S_x$ changes only finitely many time in each bounded interval of time $[0, T]$ and only a single site changes each time.

3. *Zero-temperature dynamics.* Every spin change either lowers the total energy or leaves it unchanged.

**Theorem 2.** *Suppose that the stochastic process $\hat{S}$ and (nearest-neighbor) Hamiltonian on $\mathbb{Z}^d$ satisfy hypotheses 1,2 and 3 above. If in addition either*

4a. $E(\max_{\{S_z : ||z-x|| \leq 1\}} |V_x(\{S_z : ||z - x|| \leq 1\})|) < \infty$

*or*

4b. *the strong percolation hypothesis is valid,*

*then almost surely every site has only finitely many changes that strictly lower the energy.*

The proof of Theorem 2 is given in Section 2 below.

**Remark 1.5.** *The nearest-neighbor assumption on $\mathcal{H}$ can be completely relaxed under a finite mean energy density hypothesis like 4a. Without finite mean energy density, one can also certainly extend the percolation hypothesis to cover finite range interactions, at the cost of some complications. It also appears that one can do a percolation Lyapunov function approach beyond the case of strictly finite range, but we have not investigated that thoroughly.*

**Remark 1.6.** *The proof of Theorem 2 makes clear that the part of the percolation hypothesis that $\lambda(K) \to \infty$ as $K \to \infty$ can be weakened to $\lim_{K \to \infty} \lambda(K) > (2d+1) \ln |\mathcal{S}_0|$.*

**Remark 1.7.** *We can weaken hypothesis 1 of invariance under all $\mathbb{Z}^d$-translations to invariance under a subgroup of $\mathbb{Z}^d$, such as $k_1 \mathbb{Z} \times \cdots \times k_d \mathbb{Z}$. In a natural way, by using $k_1 \times k_2 \times \cdots \times k_d$ block variables, we can obtain a new $\mathbb{Z}^d$-lattice and a new space $\mathcal{S}_0$ such that the new $V_x$'s are translation invariant. We stress that the finite mean energy or percolation hypothesis must be verified for the new $V_x$'s.*

We conclude with a few examples to show that the extensions of Theorem 2 beyond Theorem 1 allow one to deal with many interesting models.

A. *More than two states.* Potts models have $\mathcal{S}_0 = \{1, \ldots, Q\}$ and replace the Ising interaction form $S_x S_y$ of (2) with a Kronecker delta interaction $\delta_{S_x, S_y}$. Continuous-time Markovian dynamics may be described in terms of rate one Poisson "update" clocks at each site $x$ and replacing $S_x = \eta$ with $S_x = \eta'$ with probability $p(\eta, \eta')$ when the clock rings. In these models, the zero-temperature dynamics used is often not the one of Theorem 1, where $p(\eta, \eta') = 1$ or $1/2$ or $0$ according to whether $\Delta_x \mathcal{H}$ is $< 0$ or $= 0$ or $> 0$, but rather (see, e.g., [5]) one corresponding to uniform choice from among those $\eta'$ (including $\eta' = \eta$) that minimize $\Delta_x \mathcal{H}$. Of course Theorem 2 only requires that $p(\eta, \eta') = 0$ if $\Delta_x \mathcal{H} > 0$.

B. *Zero-temperature non-Markovian processes.* We modify the previous Potts models by changing the rule to update a site. Every time $t$ that a site is updated it takes a value $\eta'$ from all those that do not increase the energy, with probability depending on the fractions of time in $[0, t]$ spent in those values. This process is evidently not Markovian in general because the update depends on the history.

C. *Non-i.i.d. disorder.* In statistical mechanics there are many models in which random couplings are not independent. For example the (nearest-neighbor) Hopfield model on $\mathbb{Z}^d$ (see, e.g., [2]) with $M$ patterns has

$$J_{\{x,y\}} = \sum_{i=1}^{M} \xi_x^{(i)} \xi_y^{(i)},$$

where the $\xi_x^{(i)}$'s are i.i.d. $\pm 1$ random variables with a symmetric Bernoulli distribution. We also note that the $V_x$'s of (9) need not be sums of two-body interactions of the form $V_{\{x,z\}}(S_x, S_z)$ but can have quite general multi-body dependence.

## 2 Proofs

We begin by proving Theorem 1 as a corollary of Theorem 2 since a direct proof of Theorem 1 does not provide much simplification. Then we prove Theorem 2 in two parts corresponding respectively to hypotheses 4a and 4b.

**Proof of Theorem 1.** Hypotheses 1, 2 and 3 of Theorem 2 are obviously valid under the assumptions of Theorem 1 and thus it suffices to show that the strong percolation hypothesis 4b is also valid. Since $K^*_{\{x,y\}} = |J_{\{x,y\}}|$, we see that the percolation model of 4b is simply independent bond percolation of density $p = P(|J_{\{x,y\}}| > K)$. Since $p \to 0$ as $K \to \infty$, the large deviation bound (8) follows from standard percolation results (see, e.g., Sec 6.3 of [8]) and it is fairly easy to show that $\lambda(K) \to \infty$ as $K \to \infty$ ($p \to 0$) by, e.g., lattice animal considerations (see, e.g., Sec. 4.2 of [8]). □

**Proof of Theorem 2 (with finite mean energy density).** In this case we can essentially copy the corresponding arguments in [9] (see Theorem 3 there). We present them here since the proof under the strong percolation hypothesis is an extension of these

arguments. Let

$$\mathcal{E}(t) = E(V_x(S^t)), \tag{10}$$

which, by translation-invariance, is independent of $x$, and by the finite mean energy density hypothesis 4a is in some bounded interval $[-\mathcal{E}_0, \mathcal{E}_0]$ for all $t \geq 0$. Let $\Delta_x^t V_y$ denote the sum of all changes in $V_y(S^{t'})$ caused by flips of $S_x^{t'}$ (from any value $\eta$ to any other value $\eta'$ in $\mathcal{S}_0$) for $t' \in (0, t)$. Thus the total energy change caused by such flips is

$$\Delta_x^t \mathcal{H} = \sum_{y: ||y-x|| \leq 1} \Delta_x^t V_y; \tag{11}$$

by the zero-temperature hypothesis 3, $\Delta_x^t \mathcal{H} \leq 0$. We denote by $\mathcal{N}_x^t(\varepsilon)$ the number of those flips for which $\Delta_x \mathcal{H} \leq -\varepsilon$.

We will show that for any $\varepsilon > 0$,

$$\mathcal{E}(t) - \mathcal{E}(0) = E(\Delta_x^t \mathcal{H}) \leq -\varepsilon E(\mathcal{N}_x^t(\varepsilon)). \tag{12}$$

The inequality is obvious while the equality can be derived as follows (by means of the "Mass Transport Principle" as in [3]): Utilizing (10) and

$$V_x(S^t) - V_x(S^0) = \sum_{y: ||x-y|| \leq 1} \Delta_y^t V_x, \tag{13}$$

taking expectations and using translation invariance with the origin of $\mathbb{Z}^d$ denoted $O$, and finally using (11), we have

$$\mathcal{E}(t) - \mathcal{E}(0) = \sum_{y: ||x-y|| \leq 1} E(\Delta_y^t V_x) = \sum_{y: ||x-y|| \leq 1} E(\Delta_O^t V_{x-y})$$

$$= \sum_{y: ||x-y|| \leq 1} E(\Delta_O^t V_{y-x}) = \sum_{y: ||x-y|| \leq 1} E(\Delta_x^t V_y) \tag{14}$$

$$= E(\Delta_x^t \mathcal{H}).$$

Since $\Delta_x^t \mathcal{H} \leq 0$ and $\mathcal{E}(t) - \mathcal{E}(0) \geq -2\mathcal{E}_0$, we may let $t \to \infty$ in (12) to conclude that

$$E(\mathcal{N}_x^\infty(\varepsilon)) \leq \varepsilon^{-1} |\mathcal{E}(\infty) - \mathcal{E}(0)| \leq 2\varepsilon^{-1} \mathcal{E}_0, \tag{15}$$

so that for every $x$, and every $\varepsilon > 0$, a.s. $\mathcal{N}_x^\infty(\varepsilon) < \infty$. The conclusion of Theorem 2 is that a.s. $\mathcal{N}_x^\infty(0+) < \infty$. To see that this follows from $\mathcal{N}_x^\infty(\varepsilon) < \infty$ for every (deterministic) $\varepsilon > 0$, note that $\mathcal{N}_x^\infty(0+) = \mathcal{N}_x^\infty(\bar{\varepsilon}_x)$ where $\bar{\varepsilon}_x$ is the *random variable*,

$$\bar{\varepsilon}_x = \min\{\phi(\eta, \eta', (S_y : ||y-x|| = 1)) : \eta, \eta', S_y \text{ are } \pm 1 \text{ and } \phi > 0\} \tag{16}$$

with

$$\phi(\eta, \eta', (S_y : ||y - x|| = 1)) = \sum_{y:\,||y-x||\leq 1} (V_y(S_x = \eta) - V_y(S_x = \eta')). \qquad (17)$$

Since $\bar{\varepsilon}_x$ is the minimum of finitely many strictly positive random variables (and in case $\phi = 0$ above for *all* choices of $\eta$, $\eta'$, $S_y$, we set $\bar{\varepsilon}_x = 1$), we have $P(\bar{\varepsilon}_x > 0) = 1$ which completes the proof. □

**Proof of Theorem 2 (with strong percolation hypothesis valid).** In the absence of finite mean energy density, our strategy is to replace $V_x(S^t)$ in (10) with a different and percolation-theoretic function $\hat{\mathcal{L}}_x(S^t)$, defined in a translation invariant way, so that $L(t) = E(\hat{\mathcal{L}}_x(S^t))$ can serve as a Lyapunov function density in place of the mean energy density $\mathcal{E}(t)$. Like in the finite mean energy density proof, for this to work, it suffices for $L(t)$ to satisfy:

(i) $L(t) \in [-L_0, L_0]$ for all $t \geq 0$ with $0 < L_0 < \infty$, and

(ii) $L(t)$ is nonincreasing and $L(t) - L(0) \leq -\varepsilon E(\mathcal{N}_x^t(\varepsilon))$ for all small $\varepsilon > 0$ and arbitrary $t \geq 0$.

Here $\mathcal{N}_x^t(\varepsilon)$ is defined, as before, in terms of energy changes $\Delta_x \mathcal{H}$ (and not in terms of Lyapunov function changes).

The Lyapunov density $\hat{\mathcal{L}}_x(S)$, and the (formal) Lyapunov function $\mathcal{L} = \sum_{x \in \mathbb{Z}^d} \hat{\mathcal{L}}_x$, will *not* have only nearest-neighbor dependence but in fact will depend on $\{S_z : z \in \bar{C}_x^K\}$ for an appropriately chosen $K$, where $\bar{C}_x^K$ is the "closure" of the percolation cluster $C_x^K$ that contains $x$. The percolation clusters $C$ for a given $K$ were defined in Sec. 1 as part of the strong percolation hypothesis and the closure $\bar{C}$ of $C$ is defined as $\{y : ||y - x|| \leq 1 \text{ for some } x \in C\}$. To obtain (i) and (ii), we claim that it suffices if

(i') $$E(\max_S |\hat{\mathcal{L}}_x(S)|) < \infty, \qquad (18)$$

and for any flip of $S_x$ with $\Delta_x \mathcal{H} \leq 0$, the corresponding change $\Delta_x \mathcal{L}$ of the Lyapunov function $\mathcal{L}$ satisfies, for some deterministic $\bar{K} > 0$,

(ii') $$\Delta_x \mathcal{L} \leq \begin{cases} \Delta_x \mathcal{H}, & \text{if } -\bar{K} \leq \Delta_x \mathcal{H} \leq 0, \\ -\bar{K}, & \text{if } \Delta_x \mathcal{H} \leq -\bar{K}. \end{cases} \qquad (19)$$

Condition (i) is then immediate from (i') while (ii) follows from the extension of (12) to $\Delta_x^t \mathcal{L}$, the total change in $\mathcal{L}$ caused by flips of $S_x^{t'}$ for $t' \in (0, t)$:

$$L(t) - L(0) = E(\Delta_x^t \mathcal{L}) \leq -\varepsilon E(\mathcal{N}_x^t(\varepsilon)) \text{ for } \varepsilon \leq \bar{K}. \qquad (20)$$

In (20), the inequality is a consequence of (ii') and the equality follows from translation invariance, exactly as for $\mathcal{E}(t) - \mathcal{E}(0)$ (see (14)). To complete the proof, we need to construct $\hat{\mathcal{L}}_x$ (in a translation invariant way) so that (i') and (ii') are valid.

To begin that construction, we choose $K$ sufficiently large such that all clusters $C_y^K$ are finite a.s.; this can be done by the strong percolation hypothesis. Then we may rewrite the Hamiltonian (9) as

$$\mathcal{H}(S) = \sum_{y \in \mathbb{Z}^d} |C_y^K|^{-1} V_{C_y^K}(S) = \sum_{C \in \mathcal{C}^K} V_C(S), \tag{21}$$

where

$$V_C(S) = \sum_{z \in C} V_z(S) \tag{22}$$

and $\mathcal{C}^K$ denote the collection of all clusters.

Our Lyapunov function density $\hat{\mathcal{L}}_y$ will be of the form

$$\hat{\mathcal{L}}_y = |C_y^K|^{-1} \tilde{\mathcal{L}}_{C_y^K}, \tag{23}$$

where $\tilde{\mathcal{L}}_C$ is closely related to $V_C$, but with modifications made to yield both (i') and (ii'). Like $V_C$, $\tilde{\mathcal{L}}_C(S)$ depends only on $\{S_z : z \in \bar{C}\}$ and thus is a function on the finite space $\mathcal{S}_0^{\bar{C}}$. Let $S_{(1)}, \ldots, S_{(M_C)}$, with $M_C = |\mathcal{S}_0|^{|\bar{C}|}$, be the spin configurations in this finite space ordered so that $\mathcal{V}_{(i)} \equiv V_C(S_{(i)})$ has $\mathcal{V}_{(1)} \leq \mathcal{V}_{(2)} \leq \cdots \leq \mathcal{V}_{(M_C)}$. Then we define $\tilde{\mathcal{L}}_C$ on this space by setting $\tilde{\mathcal{L}}_C(S_{(1)}) = 0$ and

$$\tilde{\mathcal{L}}_C(S_{(i+1)}) = \sum_{j=1}^{i} \min(\mathcal{V}_{(j+1)} - \mathcal{V}_{(j)}, 4dK). \tag{24}$$

This definition gives an $\tilde{\mathcal{L}}_C$ that is a strictly increasing function of $V_C(S)$ in such a way that the following two properties are valid:

(i'') $$0 \leq \tilde{\mathcal{L}}_C \leq 4dK|\mathcal{S}_0|^{|\bar{C}|} \tag{25}$$

and, for any $S$, $S'$, the difference $\Delta\tilde{\mathcal{L}}_C \equiv \tilde{\mathcal{L}}_C(S') - \tilde{\mathcal{L}}_C(S)$ relates to the corresponding change $\Delta V_C$ according to:

(ii'') $$\Delta\tilde{\mathcal{L}}_C \begin{cases} = \Delta V_C, & \text{if } |\Delta V_C| \leq 4dK, \\ \geq 4dK, & \text{if } \Delta V_C \geq 4dK, \\ \leq -4dK, & \text{if } \Delta V_C \leq -4dK. \end{cases} \tag{26}$$

Now

$$\mathcal{L} = \sum_y \hat{\mathcal{L}}_y = \sum_y |C_y^K|^{-1} \tilde{\mathcal{L}}_{C_y^K} = \sum_{C \in \mathcal{C}} \tilde{\mathcal{L}}_C, \tag{27}$$

so that

$$\Delta_x \mathcal{L} = \Delta_x \tilde{\mathcal{L}}_x + \Delta_x \mathcal{L}_x^*, \tag{28}$$

where $\tilde{\mathcal{L}}_x = \tilde{\mathcal{L}}_{C_x^K}$ and $\mathcal{L}_x^*$ is the sum of a (random) number, between 0 and $2d$, of distinct $\tilde{\mathcal{L}}_C$'s corresponding to clusters $C \in \mathcal{C}^K$ with $x \in \bar{C} \setminus C$. For each such $C$ and each $z \in C$, either $\Delta_x V_z = 0$ (if $||z - x|| > 1$) or $|\Delta_x V_z| \leq K$ (if $||z - x|| = 1$) by (6) and the fact that $x$ and $z$ are in different percolation clusters so $K_{x,z} \leq K^*_{\{x,z\}} \leq K$. Since $x$ has exactly $2d$ nearest neighbors so there are at most $2d$ sites $z$ other than $x$ with $|\Delta_x V_z| \neq 0$, we see that $|\Delta_x V_C| \leq 2dK$ and hence, by (26), $\Delta_x \tilde{\mathcal{L}}_C = \Delta_x V_C$. Summing over such $C$'s, we get

$$\Delta_x(\mathcal{L}_x^*) = \Delta_x \sum_{y \notin C_x^K} V_y \in [-2dK, 2dK]. \tag{29}$$

If $\Delta_x \tilde{\mathcal{L}}_x = \Delta_x V_{C_x^K}$, then $\Delta_x \mathcal{L} = \Delta_x(\tilde{\mathcal{L}}_x + \tilde{\mathcal{L}}_x^*) = \Delta_x \mathcal{H}$ ($\leq 0$) and (ii′) is valid (for any choice of $\bar{K} > 0$). If $\Delta_x \tilde{\mathcal{L}}_x \neq \Delta_x V_{C_x^K}$, then by (ii″) (and (29) and the fact that $\Delta_x \mathcal{H} \leq 0$) we must have $\Delta_x V_{C_x^K} \leq -4dK$ and $\Delta_x \tilde{\mathcal{L}}_x \leq -4dK$ so that by (29)

$$\Delta_x \mathcal{H}, \ \Delta_x \mathcal{L} \leq -4dK + 2dK = -2dK; \tag{30}$$

thus (ii′) is valid with $\bar{K} = 2dK$.

It remains to show that for some choice of (large) $K$, the finite mean condition (i′) will be valid. Here is where we use the strong percolation hypothesis. By (23) and (i″),

$$E(\max_S |\hat{\mathcal{L}}_x(S)|) \leq 4dK E(|\mathcal{S}_0|^{|\bar{C}_y^K|}). \tag{31}$$

Using the crude bound that $|\bar{C}_y^K \setminus C_y^K| \leq 2d|C_y^K|$, we have from (31) that

$$E(\max_S |\hat{\mathcal{L}}_x(S)|) \leq 4dK E(e^{\alpha |C_x^K|}), \tag{32}$$

with $\alpha = (2d + 1) \ln |\mathcal{S}_0|$. By the strong percolation hypothesis the right-hand side of (32) will be finite for $K$ large enough so that the $\lambda(K)$ of the large deviation bound (8) satisfies $\lambda(K) > \alpha$. This yields (i′) and completes the proof. $\square$

**Acknowledgments.** Research partially supported by CNR (E. De Santis) and by the U.S. NSF under grants DMS-98-02310 and DMS-01-02587 (C. Newman). One of us (E. De Santis) thanks the Courant Institute for its hospitality.

# References

[1] ARRATIA, R. (1983). Site recurrence for annihilating random walks on $\mathbb{Z}_d$, *Ann. Probab.* **11** 706-713.

[2] BOVIER, A., GAYRARD, V. and PICCO, P. (1998). Typical profiles of the Kac-Hopfield model, *Mathematical Aspects of Spin Glasses and Neural Networks*, A. Bovier and P. Picco, Eds., Birkhäuser 187-241.

[3] BENJAMINI, I., LYONS, R., PERES, Y. and SCHRAMM, O. (1999). Critical percolation on any nonamenable group has no infinite clusters, *Ann. Probab.* **27** 1347-1356.

[4] CAMIA, F., DE SANTIS, E. and NEWMAN, C. M. (2001). Clusters and recurrence in the two-dimensional zero-temperature stochastic Ising model, *Ann. Appl. Probab.*, to appear.

[5] FONTES, L. R., ISOPI, M., NEWMAN, C. M. and STEIN, D. L. (2001). Aging in $1D$ discrete spin models and equivalent systems, *Phys. Rev. Lett.* **87** 110201 (1-4).

[6] FONTES, L. R., SCHONMANN, R. and SIDORAVICIUS, V. (2001). Stretched exponential fixation in stochastic Ising model at zero temperature, preprint mp-arc 01-287.

[7] GANDOLFI, A., NEWMAN, C. M. and STEIN, D. L. (2000). Zero-temperature dynamics of $\pm J$ spin glasses and related models, *Commun. Math. Phys.* **214** 373-387.

[8] GRIMMETT, G. R. (1999). *Percolation*. Second edition. *Spinger-Verlag.*

[9] NANDA, S., NEWMAN, C.M. and STEIN, D. L. (2000). Dynamics of Ising spin systems at zero temperature, in *On Dobrushin's way (from Probability Theory to Statistical Mechanics)*, R. Minlos, S. Shlosman, and Y. Suhov, Eds., *Amer. Math. Soc. Transl.* (2) **198** 183-193.